\documentclass[12pt]{article}
\usepackage{amsmath}
\usepackage{amssymb}
\usepackage{graphicx}
\usepackage[cp1251]{inputenc}

\usepackage[russian]{babel}
\newcommand{\il}[2]{\int\limits_{#1}^{#2}}
\newcommand{\ilp}[1]{\int\limits_{#1}^{+\infty}}

\newcommand{\ph}{\phantom{a}}
\newcommand{\phh}{\phantom{aaa}}

\newcommand{\sist}[2]{\left\{
\begin{array}{l}
{#1}\\
\ph\\
{#2}
\end{array}
\right.}

\begin{document}

\vskip 20pt

MSC 34C10

\vskip 10pt

\centerline{\bf Oscillation  criteria  for extended}
 \centerline{\bf linear matrix Hamiltonian systems}

\vskip 10 pt

\centerline{\bf G. A. Grigorian}
\centerline{\it Institute  of Mathematics NAS of Armenia}
\centerline{\it E -mail: mathphys2@instmath.sci.am}
\vskip 10 pt

\noindent
Abstract. The Riccati equation method is used to establish new oscillation criteria for extended linear matrix Hamiltonian systems. This method allows  to  obtain results in in a new direction, which is to break the positive definiteness condition, imposed on one of coefficients of the system.
Some examples are provided for comparing the obtained results with each other and with the result of  K. I. Al - Dosary, H. Kh. Abdullah and D. Husein.

\vskip 10 pt


\vskip 10 pt

{\bf 1. Introduction}.
Let $A(t), \ph B(t)$ and $C(t)$ be complex-valued locally integrable   matrix  functions of dimension $n \times n$ on $[t_0,\infty)$ and let Let $\mu(t)$ be real-valued locally integrable function on $[t_0,\infty), \ph I$ be the identity matric of dimension $n \times n$.
Throughout we will assume that  $B(t)$ and $C(t)$ are Hermitian,  i.e.  $B(t) = B^*(t), \ph C(t) = C^*(t), \ph t\ge t_0$, where $*$ denotes the conjugation sign. Consider the extended linear matrix Hamiltonian system
$$
\sist{\Phi' = A(t) \Phi + B(t) \Psi,}{\Psi' = C(t)\Phi + [\mu(t) I- A^*(t)] \Psi, \phh t\ge t_0.} \eqno (1.1)
$$
 By a solution $(\Phi(t), \Psi(t))$ of this system we mean a pair of absolutely continuous matrix functions $\Phi = \Phi(t)$ and $\Psi = \Psi(t)$ of dimension $n\times n$ on $[t_0,\infty)$ satisfying (1.1) almost everywhere on $[t_0, \infty)$.

{\bf Definition 1.1}. {\it A solution $(\Phi(t), \Psi(t))$ of the system (1.1) is called conjoined (or prepared, preferred), if $\Phi^*(t) \Psi(t) = \Psi^*(t) \Phi(t), \ph t\ge t_0.$}

{\bf Definition 1.2}. {\it The system (1.1) is called oscillatory if  for its every conjoined solution $(\Phi(t), \Psi(t))$  the function $\det \Phi(t)$ has arbitrarily large zeroes.}

{\bf Remark 1.1.} {\it
We claim, that if  for a solution $(\Phi(t),\Psi(t))$ of the system (1.1) the identity  $\Phi^*(t) \Psi(t) = \Psi^*(t) \Phi(t)$ holds for any  $t =t_1 \ge t_0,$ then it  holds for all $t \ge t_0$. Indeed, by (1.1) we have
$$
[\Phi^*(t) \Psi(t) - \Psi^*(t) \Phi(t)]' = \Phi^{*'}(t)\Psi(t) +\Phi^*(t) \Psi'(t) - \Psi^{*'}(t) \Phi(t) - \Psi^*(t) \Phi'(t) =
$$
$$
=[A(t) \Phi(t) + B(t) \Psi(t)]^* \Psi(t) + \Phi^*(t)[C(t)\Phi(t) + [\lambda(t) I- A^*(t)] \Psi(t)] -
$$
$$
- [C(t)\Phi(t) + [\lambda(t) I- A^*(t)] \Psi(t)]^* \Phi(t) - \Psi^*(t) [A(t) \Phi(t) + B(t) \Psi(t)], \phh t \ge t_0.
$$
After some simplifications from here we obtain
$$
[\Phi^*(t) \Psi(t) - \Psi^*(t) \Phi(t)]' = \lambda(t) [\Phi^*(t) \Psi(t) - \Psi^*(t) \Phi(t)], \phh t \ge t_0.
$$
Therefore, by virtue of the Cauchy formula we get
$$
[\Phi^*(t) \Psi(t) - \Psi^*(t) \Phi(t)]= [\Phi^*(t_1) \Psi(t_1) - \Psi^*(t_1) \Phi(t_1)] \exp\biggl\{\il{t_1}{t}\lambda(\tau)d\tau\biggr\}, \phh t \ge t_0.
$$
Hence, the claim.
}

\vskip 10pt

The oscillation property of the system (1.1) has been studied in [1] for the case of its real-valued continuous coefficients. It was established in [1] that the system is reducible to an Hamiltonian system by an appropriate linear transformation. By this way it was obtained the following result. Let $A(t), \ph B(t), \ph C(t)$ and $\mu(t)$ be real-valued continuous on $[t_0,\infty)$. Denote $A_1 = A- BK, \ph B_1 = \alpha e^{\il{\ph}{\ph} \mu} B, \ph C_1 = \alpha e^{-\il{\ph}{\ph} \mu}(KA -KBK + C -\mu K +A^*K),$ where $\alpha = \sist{1, \ph if  \ph B > 0,}{-1, \ph if \ph B< 0,} \ph K$ is an constant nonzero symmetric $n\times n$ matrix.

{\bf Theorem 1.1 ([1, Proposition 4])}. {\it Let $A(t), \ph B(t)$ and $C(t)$ be real-valued continuous functions on $[t_0,+\infty)$ and  $B(t)$ be positive definite for all $t \ge t_0$. If there exists a positive linear functional $g$ on the space of matrices of dimension $n\times n$ such that
$$
\lim\limits_{t \to \infty}\il{t_0}{t}\frac{\alpha e^{\il{}{}\mu}d s}{g[B^{-1}(s)]} = \infty,
$$
$$
\lim\limits_{t\to +\infty}g\left[ -\il{t_0}{t} \Bigl(C_1(s) + A_1^*(s) B_1^{-1}(s) A_1(s)\Bigr) ds - B_1^{-1}(t) A_1(t)\right] = \infty,
$$
then the system (1.1) is oscillatory.
}

This theorem generalizes a result of work [10] (see [10, Theorem 2.1]). Note that in [1] unlike of the established above (Remark 1.1) the existence of conjoined solutions of the system it was established by indirect way, i. e. it was established that a solution of the system (1.1) is conjoined if and only if the corresponding solution to the transformed Hamiltonian system is conjoined (see [1, p. 26, Proposition 2]).
The oscillation problem for linear matrix Hamiltonian systems is that of finding explicit conditions on the coefficients of the system, providing its oscillation. This is an important problem
of qualitative theory of differential equations and  is a subject of many investigations  (see e.g., [1, 2, 4, 5, 12-17] and cited works therein).
In this paper we use the Riccati equation method for obtaining some new oscillation criteria for the system (1.1). We show that one of the obtained results is an extension of Theorem 1.1.
 On examples we compare the obtained oscillation criteria with each other and with Theorem 1.1 and demonstrate their applicability.

{\bf 2. Auxiliary propositions.}
Let $a_{jk}(t), \ph j,k =1,2$ be real-valued locally integrable functions on $[t_0,\infty)$. Consider the linear system of ordinary differential equations
$$
\sist{\phi' = a_{11}(t) \phi + a_{12}(t) \psi,}{\psi' = a_{21}(t) \phi + a_{22}(t) \psi, \ph t\ge t_0.} \eqno (2.1)
$$
and the corresponding Riccati equation
$$
y' + a_{12}(t) y^2 + E(t) y -a_{21}(t) = 0, \phh t \ge t_0, \eqno (2.2)
$$
where $E(t) \equiv a_{11}(t) - a_{22}(t), \ph t \ge t_0$. By a solution of the system (2.1) we mean an ordered pair $(\phi(t), \psi(t))$ of absolutely continuous functions $\phi(t), \ph \psi(t)$ on $[t_0,\infty)$, satisfying for $\phi =\phi(t), \ph \psi = \psi(t)$  the system (2.1) almost everywhere on $[t_0,\infty)$.
All solutions $y(t)$ of Eq. (2.2), existing on any interval $[t_1,t_2) \subset [t_0,\infty)$ are connected with solutions $(\phi(t),\psi(t))$  of the system (2.1) by relations (see [6], [7])
$$
\phi(t) = \phi(t_1)\exp\biggl\{\il{t_1}{t}[a_{12}(\tau) y(\tau) + a_{11}(\tau)]d\tau\biggr\}, \ph \phi(t_1) \ne 0, \phh  \psi(t) = y(t) \phi(t),  \eqno (2.3)
$$
$t \in [t_1,t_2).$

{\bf Definition 2.1.} {\it The system (2.1) is called oscillatory if for its  every solution  $(\phi(t), \psi(t))$ the function $\phi(t)$ has arbitrarily large zeroes.}

\vskip 5pt

{\bf Remark 2.3}. {\it Explicit oscillatory criteria for the system  (2.1) are obtained in [6].}

{\bf Theorem 2.1.} {\it Let the following conditions be satisfied.

\noindent
$a(t) \ge 0, \ph t \ge t_0.$
$$
\ilp{t_0}a_{12}(t)a_{12}(t)\exp\biggl\{-\il{a}{t}E(\tau) d \tau\biggr\} d t = - \ilp{t_0}a_{21}(t)\exp\biggl\{\il{a}{t}E(\tau) d \tau\biggr\}d t = \infty.
$$
Then the system (2.1) is oscillatory.
}

Proof. By analogy with the proof of Corollary 3.1 from [9] (see also [5, Theorem 2.4]).

{\bf Definition 2.2.} {\it An interval $[t_1,t_2)\subset [t_0,\infty)$ is called the maximum existence interval for a solution $y(t)$ of Eq. (2.2)
if $y(t)$ exists on $[t_1,t_2)$ and cannot be continued to the right from $t_2$ as a solution of that equation.
}

{\bf Lemma 2.1.} {\it Let $y(t)$ be a solution of Eq. (2.2) on $[t_1,t_2) \subset [t_0,+\infty)$, and let $t_2< +\infty$. If  $a(t) \ge 0, \ph t \in [t_1,t_2)$ and  the function $F(t) \equiv \il{t_1}{t}a_{12}(\tau) y(\tau) d\tau, \ph t \in [t_1,t_2)$ is bounded from below on $[t_1,t_2)$, then $[t_1,t_2)$ is not the maximum existence interval for $y(t)$.
}

Proof. By analogy with the proof of Lemma 2.1 from [8].

Let $e(t)$ and $e_1(t)$ be real-valued functions on $[t_0,\infty)$ and let $e(t)$ be locally integrable and $e_1(t)$ be absolutely continuous on $[t_0,\infty)$. Consider the Riccati integral equations
$$
y(t) + \il{t_0}{t} a(\tau) y^2(\tau) d\tau + e(t) = 0, \phh t \ge t_0, \eqno (2.4)
$$
$$
y(t) + \il{t_0}{t} a(\tau) y^2(\tau) d\tau + e_1(t) = 0, \phh t \ge t_0, \eqno (2.5)
$$

{\bf Lemma 2.2.} {\it Let $y_0(t)$ be a solution of Eq. (2.4) on $[t_0,t_1)$. If $a(t) \ge 0, \ph e(t) > e_1(t) > 0, \ph t\in[t_0,t_1)$, then Eq. (2.5) has a solution $y_1(t)$  on $[t_0,t_1)$ and
$$
y_1(t) > y_0(t), \phh t\in [t_0,t_1). \eqno (2.6)
$$
}

Proof. Since $a(t) \ge 0, \ph e(t) > 0, \ph t\in [t_0,t_1)$ by (2.4)
$$
y_0(t) < 0, \phh t \in [t_0,t_1). \eqno (2.7)
$$
Let $y_1(t)$ be a solution of Eq. (2.5) and let $[t_0,t_2)$ be its maximum existence interval. Show that
$$
t_2 \ge t_1. \eqno (2.8)
$$
Suppose
$$
t_2 < t_1. \eqno (2.9)
$$
Show that
$$
y_1(t) > y_0(t), \ph t\in [t_0,t_2). \eqno (2.10)
$$
Suppose that this is false. By (2.4) and (2.5) from the conditions $e(t) > e_1(t) > 0, \linebreak t\in [t_0,t_1)$ of the lemma it follows that $y_1(t_0) = - e_1(t_0) > - e(t_0) = y_0(t_0)$. Then there exists $t_3\in(t_0,t_2)$ such that
$$
y_1(t) > y_0(t), \ph t_0 \le t < t_3, \eqno (2.11)
$$
$$
y_1(t_3) = y_0(t_3). \eqno (2.12)
$$
On the other hand by (2.4) and (2.5) we have $y_1(t_3) - y_0(t_3) = \il{t_0}{t_3}a(\tau)[y_0^2(\tau) - y_1^2(\tau)]d\tau + e(t_3) - e_1(t_3)$. This together with (2.7), (2.11)  and the condition $a(t) \ge 0, \ph t \in [t_0,t_1)$  of the lemma implies that $y_1(t_3) > y_0(t_1)$, which contradicts (2.11). The obtained contra-\linebreak diction proves (2.10). Obviously $y_1(t)$ is a solution of the Riccati equation (recall that $e_1(t)$ is absolutely continuous)
$$
y' + a(t) y^2 + e'_1(t) = 0, \phh t \ge t_0
$$
on $[t_0,t_1)$. Then by Lemma 2.1  and (2.3) from the condition $a(t)\ge 0, \ph t\in [t_0,t_1)$ of the lemma and from (2.2) it follows that $[t_0,t_2)$ is not the maximum existence interval for $y_1(t)$, which contradicts our supposition. The obtained contradiction proves (2.8). From (2.8) and (2.10) it follows existence $y_1(t)$ on $[t_0,t_1)$ and the inequality (2.6). The lemma is proved.

{\bf Lemma 2.3.} {\it  For any two square matrices  $M_k \equiv (m_{ij}^l)_{i,j=1}^n, \ph l=1,2$  the equality
$$
tr(M_1 M_1) = tr (M_2 M_1)
$$
is valid.}

Proof.  We have $tr (M_1 M_2) = \sum\limits_{j=1}^n(\sum\limits_{k=1}^n m_{jk}^1 m_{kj}^2) = \sum\limits_{k=1}^n(\sum\limits_{j=1}^n m_{jk}^1 m_{kj}^2) = \sum\limits_{k=1}^n(\sum\limits_{j=1}^n m_{kj}^2 m_{jk}^1) = tr (M_2 M_1).$ The lemma is proved.

Let $p(t)$ be positive  continuously differentiable function on $[t_0,\infty)$.
The substitution
$$
\Psi = p(t)Y\Phi. \eqno (2.13)
$$
in the system (1.1) yields
$$
\sist{\Phi' = \Bigl[A(t) + p(t) B(t) Y\Bigr]\Phi,}{\Bigl[Y'+ p(t)Y B(t) Y + A^*(t) Y + Y \bigl[A(t) + \frac{p'(t)}{p(t)}I + \mu(t)I\bigr] - \frac{1}{p(t)}C(t)\Bigr]\Phi = 0, \ph t \ge t_0.}
$$
It follows from here and from (2.13) that all solutions $Y(t)$ of the matrix Riccati equation
$$
Y'+ p(t)Y B(t) Y + A^*(t) Y + Y \Bigl[A(t) + \frac{p'(t)}{p(t)}I + \mu(t)I\Bigr] - \frac{1}{p(t)}C(t)=0, \ph t \ge t_0, \eqno (2.14)
$$
existing on any interval $[t_1,t_2) \subset [t_0,\infty)$, are connected with solutions $(\Phi(t),\Psi(t))$ of the system (1.1) by relations.
$$
\Phi'(t) = [A(t) + p(t) B(t) Y(t)] \Phi(t), \phh \Psi(t) = p(t) Y(t) \Phi(t), \ph t \in [t_1,t_2). \eqno (2.15)(2.19)
$$
(by a solution of Eq. (2.14)  on $[t_1,t_2) \subset [t_0,\infty)$ we mean an absolutely continuous matrix function on $[t_1,t_2)$, satisfying (2.14) almost everywhere on $[t_1,t_2)$).
Let us set $A_1(t) \equiv A(t) + \frac{1}{2}\bigl(\frac{p'(t)}{p(t)} + \mu(t)\bigr)I, \ph B_1(t)\equiv p(t)B(t), \ph C_1(t) \equiv \frac{1}{P(t)} C(t).$ Then Eq. (2.14) can be rewritten in the form
$$
Y'+ Y B_1(t) Y + A_1^*(t) Y + Y A_1(t) - C_1(t)=0, \ph t \ge t_0. \eqno (2.16)
$$
For any matrix $M$ of dimension $n\times n$ denote by $\lambda_1(M), \dots, \lambda_n(M)$ the eigenvalues of $M$, and if they are real-valued then we will assume that they are ordered by
$$
\lambda_1(M) \le \dots,\le \lambda_n(M).
$$
The nonnegative (positive) definiteness of any Hermitian matrix will be denoted by $H\ge~ 0  \linebreak (> 0)$. By $I$ it will be denoted the identity matrix of dimension $n\times n$.

{\bf Lemma 2.4.} {\it For a matrix $S$ of dimension $n\times n$ and any Hermitian matrix $H \ge 0$ of the same dimension the inequality
$$
tr(SHS^*) \ge \frac{\lambda_1(H)}{n}\biggl\{\biggl[tr \biggl(\frac{S + S^*}{2}\biggr)\biggr]^2 + \biggl[tr \biggl(\frac{S - S^*}{2 i}\biggr)\biggr]^2\biggr\}
$$
is valid.
}

Proof. Let $U_H$ be a unitary matrix such that
$$
U_h H U^*_H = diag \{\lambda_1(H),\dots,\Lambda_n(H)\} \stackrel{def}{=}  H_0.
$$
and let
$$
S_H = U_H S U^*_H = (s_{jk})_{j,k=1}^n.
$$
Then
$$
tr (S H S^*) = tr (S_H H_0 S_H) = \sum\limits_{j,k=1}^n \lambda_k(H) s_{jk} \overline{s_{jk}}. \eqno (2.17)
$$
Since $H\ge 0$ we have $\lambda_n(H)\ge\dots \lambda_1(H) \ge 0$. This together with (2.17) implies
$$
tr (S H S^*) \ge \lambda_1(H) \sum\limits_{j,k=1}^ns_{jk} \overline{s_{jk}} \ge \lambda_1(H)\sum\limits_{j=1}^ns_{jj} \overline{s_{jj}} =
$$
$$
 =\lambda_1(H)\biggl[\sum\limits_{j=1}^n\biggl(\frac{s_{jj} + \overline{s_{jj}}}{2}\biggr)^2 + \sum\limits_{j=1}^n\biggl(\frac{s_{jj} - \overline{s_{jj}}}{2 i}\biggr)^2\biggr] \ge
$$
$$
\ge\frac{\lambda_1(H)}{n}\biggl\{\biggl[\sum\limits_{j=1}^n\frac{s_{jj} + \overline{s_{jj}}}{2}\biggr]^2 + \sum\limits_{j=1}^n\biggl\{\biggl[\frac{s_{jj} - \overline{s_{jj}}}{2 i}\biggr]^2\biggr\} =
$$
$$
=\frac{\lambda_1(H)}{n}\biggl\{\biggl[tr \biggl(\frac{S + S^*}{2}\biggr)\biggr]^2 + \biggl[tr \biggl(\frac{S - S^*}{2 i}\biggr)\biggr]^2\biggr\}.
$$
The lemma is proved.

Let $g$ be a positive linear functional on the space of matrices of dimension $n\times n$. For any matrix $M\ge 0$ of dimension $n\times n$ set
$$
\nu_g(M)\equiv\sist{0, \ph if \ph \det M = 0,}{\{g(M^{-1})\}^{-1}, \ph if \ph \det M \ne 0.}
$$

It is known that for every Hermitian matrix $D \ge 0$ of dimension $n\times n$ the estimates
$$
\lambda_1(D) \le g(D) \le \lambda_n(D) \eqno (2.18)
$$
are valid for every positive linear functional $g$ (see [17]). Then from the relation $g(B^{-1}(t) \ge \lambda_1(B(t)), \ph B(t) > 0$ it follows that
$$
\nu_g(B(t)) \le \lambda_1(B(t)) \le tr B(t), \phh t \ge t_0. \eqno (2.19)
$$
provided $B(t) \ge 0, \ph t \ge t_0.$ Hence, $\nu_g(B(t)), \ph t\ge t_0$ is always locally integrable for $B(t)\ge 0, \ph t \ge t_0$.

\vskip 10pt

{\bf 3. Oscillation criteria.} Hereafter by the satisfiability of a relation $\mathcal{P}$ (equality, inequality) on any interval we will mean (if it is necessary) the satisfiability of  $\mathcal{P}$ almost everywhere on that interval.
For any matrix function $P(t), \ph t\ge t_0$ of dimension $n\times n$ set
$$
J_P(t) \equiv -\il{t_0}{t} [C(\tau) + A^*(\tau) P(\tau)]d \tau - P(t), \ph t \ge t_0.
$$
Denote by $\mathcal{M}_\mathbb{R}$ the set of matrices of dimension $n\times n$ with real entries.
Denote
$$
J(t) \equiv \frac{1}{p(t)} tr\biggl[ \biggl(\frac{A(t) + A^*(t)}{2} + \Bigl(\frac{p'(t)}{p(t)} + \mu(t)\Bigr)I\biggr)B^{-1}(t)\biggr] -
$$
$$
- \il{t_0}{t}\frac{1}{p(t)}tr \Biggl\{A(\tau) B^{-1}(\tau) A^*(\tau) +\frac{1}{2}\Bigl(\frac{p'(\tau)}{p(\tau)}\Bigr)\Bigl[B^{-1}(\tau)A^*(\tau) + A(\tau)B^{-1}(\tau)\Bigr] +
$$
$$
+ \frac{1}{4p(\tau)}\Bigl(\frac{p'(\tau)}{p(\tau)} + \mu(\tau)\Bigr)^2 B^{-1}(\tau) + \frac{1}{p(\tau)}C(\tau)\Biggr\} d \tau
 +\il{t_0}{t}\frac{p(\tau)\lambda_1(B(\tau))}{n}\biggl[ tr\biggl(\frac{A(\tau) - A^*(\tau)}{2 i}\biggr)\biggr]^2 d \tau,
$$
$t \ge t_0.$.

{\bf Theorem 3.1.} {\it Let the the functions  $tr \Bigl[A(t) B^{-1}(t) A^*(t)\Bigr], \ph \lambda_1(B(t))\biggl[ tr \bigl(A(t) - A^*(t)\bigr)\biggr]^2, \linebreak   t \ge t_0$ be locally integrable and let the following conditions be satisfied.

\noindent
I') $B(t) > 0, \ph t \ge t_0$.

\noindent
IV) $\ilp{t_0}p(\tau) \lambda_1(B(\tau)) d \tau = \lim\limits_{t \to +\infty} J(t) = \infty.$

\noindent
Then the system (1.1) is oscillatory.
}

Proof. Suppose the system (1.1) is not oscillatory. Then by (2.15) Eq. (2.14) has a Hermitian solution $Y(t)$ on $[t_1,\infty)$ for some $t_1 \ge t_0$. Then using I') and (2.16) we can write
$$
Y'(t) + \frac{1}{2}\Bigl\{[Y(t) + A_1(t) B_1^{-1}(t)] B_1(t) [Y(t) + B_1^{-1}(t) A_1^*(t)] +
$$
$$
 +[Y(t) + A_1^*(t) B_1^{-1}(t)] B_1(t) [Y(t) + B_1^{-1}(t) A_1(t)]\Bigr\}+
$$
$$
+\frac{A_1^*(t) - A_1(t)}{2} Y(t) + Y(t) \frac{A_1(t) - A_1^*(t)}{2} -
$$
$$
-\frac{1}{2}[A(t)_1B^{-1}(t) A_1^*(t) + A_1^*(t) B_1^{-1}(t) A_1(t)] - C_1(t) = 0, \ph t\ge t_1. \eqno (3.1)
$$
Since $Y(t)$ and $B_1^{-1}(t)$ are Hermitian by Lemma 2.4 we have
$$
tr \frac{1}{2}\Bigl\{[Y(t) + A_1(t) B_1^{-1}(t)] B_1(t) [Y(t) + B_1^{-1}(t) A_1^*(t)] +
$$
$$
 +[Y(t) + A_1^*(t) B_1^{-1}(t)] B_1(t) [Y(t) + B_1^{-1}(t) A_1(t)]\Bigr\} \ge
$$
$$
\frac{\lambda_1(B_1(t))}{2n}\biggl\{\biggl[tr\biggl(Y(t) + \frac{A_1(t) B_1^{-1}(t) + B_1^{-1}(t) A_1^*(t)}{2}\biggr)\biggr]^2 + \biggl[tr\biggl(\frac{A_1(t) B_1^{-1}(t) - B_1^{-1}(t) A_1^*(t)}{2 i}\biggr)\biggr]^2
$$
$$
+\biggl[tr\biggl(Y(t) + \frac{A_1^*(t) B_1^{-1}(t) + B_1^{-1}(t) A_1(t)}{2}\biggr)\biggr]^2 + \biggl[tr\biggl(\frac{A_1^*(t) B_1^{-1}(t) - B_1^{-1}(t) A_1(t)}{2 i}\biggr)\biggr]^2\biggr\},
$$
$t \ge t_1$. By Lemma 2.3 from here we obtain
$$
tr \frac{1}{2}\Bigl\{[Y(t) + A_1(t) B_1^{-1}(t)] B_1(t) [Y(t) + B_1^{-1}(t) A_1^(t)] +
$$
$$
 +[Y(t) + A_1^*(t) B_1^{-1}(t)] B_1(t) [Y(t) + B_1^{-1}(t) A_1(t)]\Bigr\} \ge
$$
$$
\frac{\lambda_1(B_1(t))}{n}\biggl\{\biggl[tr\biggl(Y(t) + \frac{A_1(t) +  A_1^*(t)}{2}  B_1^{-1}(t)\biggr)\biggr]^2 + \biggl[tr\biggl(\frac{A_1(t) - A_1^*(t)}{2i} B_1^{-1}(t)\biggr)\biggr]^2\biggr\}, \ph t \ge t_1.
$$
This together with (3.1) implies
$$
tr Y'(t) + \frac{\lambda_1(B_1(t))}{n}\biggl\{\biggl[tr\biggl(Y(t) + \frac{A_1(t) +  A_1^*(t)}{2}  B_1^{-1}(t)\biggr)\biggr]^2 - \phantom{aaaaaaaaaaaaaaaaaaaaaaaaaaaaa}
$$
$$
-\frac{1}{2}tr\Bigl\{[A_1(t)B_1^{-1}(t) A_1^*(t) + A_1^*(t) B_1^{-1}(t) A_1(t)] - C_1(t)\Bigr\} + \biggl[tr\biggl(\frac{A(t) - A^*(t)}{2i} B^{-1}(t)\biggr)\biggr]^2\biggr\} \le 0,
$$
$t \ge t_1$.  If we substitute $Z(t) \equiv Y(t) + \frac{A_1(t) + A_1^*(t)}{2} B_1^{-1}(t), \ph t \ge t_1$ in the above inequality  and integrate (by taking into account the condition of local integrability of the functions $tr [A(t)B^{-1}(t) A^*(t)], \ph tr \lambda_1(B(t))[tr(A(t) - A^*(t))]^2$) from $t_1$ to $t$ we obtain
$$
tr Z(t) + \il{t_1}{t}p(\tau)\frac{\lambda_1(B(\tau))}{n}[tr Z(\tau)]^2 d \tau + J(t) + c \le 0, \phh t \ge t_1, \eqno (3.2)
$$
where $c = Y(t_1) + \il{t_0}{t_1} tr [C_1(\tau) + A_1(\tau) B_1^{-1}(\tau) A_1^*(\tau)] d \tau- \il{t_0}{t_1}\biggl[tr\biggl(\frac{A_1(\tau) - A_1^*(\tau)}{2i} B_1^{-1}(\tau)\biggr)\biggr]^2d\tau = const.$
Consider the integral Riccati equations
$$
y(t) + \il{t_1}{t}p(\tau)\frac{\lambda_1(B(\tau))}{n}[y(\tau)]^2 d \tau + J(t) + c +l(t) = 0, \phh t \ge t_1, \eqno (3.3)
$$
$$
y(t) + \il{t_1}{t}p(\tau)\frac{\lambda_1(B(\tau))}{n}[y(\tau)]^2 d \tau + J(t) + c -1 = 0, \phh t \ge t_1, \eqno (3.4)
$$
 where $l(t)\equiv -tr Z(t) - \il{t_1}{t} p(\tau)\frac{\lambda_1(B(\tau))}{n}[tr Z(\tau)]^2d\tau -J(t) - c, \ph t \ge t_1$. By (3.2)
$$
l(t) \ge 0, \ph t \ge t_1. \eqno (3.5)
$$
Obviously, $y(t) =tr Z(t), \ph t \ge t_1$ is a solution of Eq. (3.3). It follows from the condition IV) of the theorem that $$
J(t) + c -1> 0, \ph t \ge t_2, \eqno (3.6)
$$
for some $t_2 \ge t_1$. Then Eq. (3.5) and (3.6)  on $[t_2,\infty)$ are equivalent to the following ones
$$
y(t) + \il{t_2}{t}p(\tau)\frac{\lambda_1(B(\tau))}{n}[y(\tau)]^2 d \tau + f(t) = 0, \phh t \ge t_2, \eqno (3.7)
$$
$$
y(t) + \il{t_2}{t}p(\tau)\frac{\lambda_1(B(\tau))}{n}[y(\tau)]^2 d \tau + f_1(t) = 0, \phh t \ge t_2, \eqno (3.8)
$$
where $f(t) \equiv \il{t_1}{t_2}p(\tau)\frac{\lambda_1(B(\tau))}{n}[tr U(\tau)]^2 d \tau  + J(t) + c +l(t), \ph f_1(t) \equiv \il{t_1}{t_2}p(\tau)\frac{\lambda_1(B(\tau))}{n}[tr U(\tau)]^2 d \tau   +  J(t) + c -1.$
Then (since $Z(t)$ is a solution of Eq. (3.7)) by Lemma 2.2 Eq. (3.8) has  a solution $y_1(t)$ on $[t_2,\infty)$. Note that $f_1(t)$ is a solution of the Riccati equation
$$
y' + P(t)\frac{\lambda_1(B(t))}{n} y^2 +f_1'(t) = 0, \ph t\ge t_0.
$$
Then by (2.7) the system
$$
\sist{\phi' = p(t)\frac{\lambda_1(B(t))}{n} \psi,}{\psi' = - f_1'(t), \ph t \ge t_2}
$$
is not oscillatory. On the other hand since $\il{t_2}{t}f_1'(\tau)d\tau = f_1(t) - f_1(t_2) \to \infty$ as $t\to \infty$ by Theorem 2.2 from the condition $I')$ it follows that the last system is oscillatory. We obtain a contradiction, completing the proof of the the theorem.

We set
$$
\nu_0(B(t)) \equiv \sist{0, \ph if \ph \det B(t) = 0,}{\frac{1}{tr (B^{-1}(t))}, \ph if \ph \det B(t) \ne 0,} \ph t \ge t_0.
$$
It is not difficult to verify that
$$
\frac{1}{tr(B^{-1}(t))} \le \lambda_1(B(t)) \le \frac{n}{tr(B^{-1}(t))}, \ph \mbox{for all} t \ge t_0, \ph \mbox{for which} \ph B(t) > 0. \eqno (3.9)
$$
Therefore Theorem 3.1 remains valid if we replace the condition $\ilp{t_0}\lambda_1(B(t)) d t = \infty$ of Theorem 3.1 by the following one $\ilp{t_0}\nu_0(B(t)) d t = \infty$. Moreover by (2.18)(2.22) in the case
$-\il{t_0}{t} \Bigl(C(s) + A^*(s) B^{-1}(s) A(s)\Bigr) ds - B^{-1}(t) A(t) \ge 0, \ph t \ge T$, for some $T\ge t_0$,
the functional $g$ in Theorem 1.1 is equivalent to the functional $tr$. Hence due to (3.9) Theorem 3.1 is a complement to Theorem 1.1.

{\bf Theorem 3.2.} {\it Let  the following conditions be satisfied.

\noindent
I') $B(t) > 0, \ph t \ge t_0.$

\noindent
V) $\ilp{t_0}\frac{p(t) d t}{tr (B^{-1}(t))} = \infty.$

\noindent
VI)  the  function   $tr [(A(t) + A^*(t)) B^{-1}(t)(A(t) + A^*(t)) + (\frac{p'(t)}{p(t)} + \mu(t))\{(A(t) + A^*(t))B^{-1}(t) + B^{-1}(t)(A(t) + A^*(t))\} +(\frac{p'(t)}{p(t)} + \mu(t))^2 I], \ph t \ge t_0$ is locally integrable on $[t_0,\infty)$ and  $\lim\limits_{t \to +\infty}- tr \biggl[\frac{2}{p(t)}(A_1(t) + A_1^*(t) +(\frac{p'(t)}{p(t)} + \mu(t)) I) B^{-1}(t) +$

\phantom{a} $+\il{t_0}{t}\frac{1}{p}(\tau)\Bigl((A(\tau) + A^*(\tau)) B^{-1}(\tau) (A(\tau) + A^*(\tau)) +  (\frac{p'(\tau)}{p(\tau)} + \mu(\tau))\{(A(\tau) + A^*(\tau))B^{-1}(\tau) + B^{-1}(\tau)(A(\tau) + A^*(\tau))\} +     4C(\tau)\Bigr) d\tau\Biggr] = \infty$.

\noindent
Then the system (1.1) is oscillatory.
}

Proof. Suppose the system (1.1) is not oscillatory. Then by (2.15) Eq. (2.14) has a solution $Y(t)$ on $[t_1,\infty)$ for some $t_1 \ge t_0$. Hence, by (2.16)
$$
Y'(t) + Y(t) B_1(t) Y(t) + A_1^*(t) Y(t) + Y(t) A_1(t) - C_1(t) = 0, \phh t \ge t_1.
$$
From here it follows
$$
tr \biggl\{Y'(t) + \Bigl[Y(t) + \frac{A_1(t) + A_1^*(t)}{2} B_1^{-1}(t)\Bigr] B_1(t) \Bigl[Y(t) +  B_1^{-1}(t) \frac{A_1(t) + A_1^*(t)}{2}\Bigr] +
$$
$$
+\frac{A_1^*(t) - A_1(t)}{2} Y(t) + Y(t) \frac{A_1(t) - A_1^*(t)}{2} - \frac{A_1^*(t) + A_1(t)}{2} B_1^{-1}(t)\frac{A_1^*(t) + A_1(t)}{2} - C_1(t)\biggr\} = 0,
$$
$t \ge t_1$. It follows from  the condition IV) that the function
 $tr [(A(t) + A^*(t)) B^{-1}(t)(A(t) + A^*(t))]$ is locally integrable over $[t_0,\infty)$.
Then substituting $Z(t) \equiv Y(t)  + \frac{A(t) + A^*(t)}{2}, \ph t \ge t_1$ in the obtained equality and integrating from $t_1$ to $t$, we obtain.
$$
tr \biggl\{Z(t)  + \il{t_1}{t}Z(\tau) B_1(\tau) Z^*(\tau) d \tau + \il{t_1}{t}\Bigl[\frac{A_1^*(\tau) - A_1(\tau)}{2} Y(\tau) + \phantom{aaaaaaaaaaaaaaaaaaaaaaaaaaaaaaaaaa}
$$
$$
\phantom{aaaaaaaaaaaaaaaaaaaaaa}+Y(\tau)\frac{A_1(\tau) - A_1^*(\tau)}{2}\Bigr] d\tau + J_1(t)\biggr\} = 0, \phh t \ge t_1, \eqno (3.10)
$$
 where
$$
J_1(t) \equiv -Y(t_1) - \frac{A_1(t) + A_1^*(t)}{2} B_1^{-1}(t) - \phantom{aaaaaaaaaaaaaaaaaaaaaaaaaaaaaaaaaaaaaaaaaaaaaa}
$$
$$
\phantom{aaaaaaaaaaaaaaaaaaaaaaaa}-\il{t_1}{t}\Bigl[\frac{A_1(\tau) + A_1^*(\tau)}{2} B_1^{-1}(\tau)\frac{A_1(\tau) + A_1^*(\tau)}{2} + C_1(\tau)\Bigr] d\tau, \ph t \ge t_1.
$$
It is not difficult to verify that
$$
J_1(t) \equiv -Y(t_1) - \frac{1}{p(t)}\biggl(\frac{A(t) + A^*(t)}{2} + \Bigl(\frac{p'(t)}{p(t)} + \mu(t)\Bigr) I\biggr) B^{-1}(t) - \phantom{aaaaaaaaaaaaaaaaaaaaaaaaaaaaaaaaaaaaaaaaaaaaaa}
$$
$$
-\il{t_1}{t}\frac{1}{p(\tau)}\Bigl[\frac{A(\tau) + A^*(\tau)}{2} B^{-1}(\tau)\frac{A(\tau) + A^*(\tau)}{2} +
$$
$$
+\frac{1}{2}\Bigl\{(A(\tau) + A^*(\tau))B^{-1}(\tau) +  B^{-1}(\tau)(A(\tau) + A^*(\tau)) + \frac{1}{4}\Bigl(\frac{p'(\tau)}{p(\tau)} + \mu(\tau)\Bigr)^2 I  +  C(\tau)\Bigr] d\tau, \ph t \ge t_1.
$$
By Lemma 2.3 we have
$$
tr \biggl[\il{t_1}{t}\Bigl[\frac{A_1^*(\tau) - A_1(\tau)}{2} Y(\tau) + Y(\tau)\frac{A_1(\tau) - A_1^*(\tau)}{2}\Bigr] d\tau\biggr] = 0, \phh t \ge t_1. \eqno (3.11)
$$
Since $Y(t), \ph A_1(t) + A_1^*(t)$ and $B_1^{-1}(t)$ are Hermitian we have also $tr (Z(t) - Z^*(t)) = 0, \linebreak t \ge~ t_1$. By Lemma 2.4 from here we obtain
$$
tr \il{t_1}{t} Z(\tau) B_1(\tau) Z^*(\tau) d \tau  \ge \il{t_1}{t}\frac{\lambda_1(B_1(\tau))}{n}\Bigl[ tr \frac{Z(\tau) + Z^*(\tau)}{2}\Bigr]^2 d \tau  = \il{t_1}{t}\frac{\lambda_1(B_1(\tau))}{n}\Bigl[ tr Z(\tau)\Bigr]^2 d \tau,
$$
$t \ge t_1$. This together with (3.10) and (3.11) implies that
$$
tr Z(t) - tr Z(t_1) + \il{t_1}{t}\frac{\lambda_1(B_1(\tau))}{n}\Bigl[ tr Z(\tau)\Bigr]^2 d \tau + tr J_1(t) \le 0, \ph t \ge t_1.
$$
Further as in the proof of Theorem 3.1 one can show that, if the conditions I'), VI) and the condition

\noindent
V') $\ilp{t_0}p(t)\lambda_1(B(t)) d t = \infty$

\noindent
are satisfied, then the system (1.1) is oscillatory. But according to (3.9) the condition V') with I') is equivalent to the condition V). Therefore under the conditions of the theorem the system (1.1) is oscillatory. The theorem is proved.

{\bf Example 3.2.} {\it Let $B(t) \equiv I, \ph A(t) \equiv A_0, \ph C(t) \equiv - A^*_0 A_0,  \ph t \ge t_0,\ph
 A_0= const$ is a real-valued matrix of dimension $n\times n, \ph p(t)\equiv 1, \ph \mu(t) \equiv 0$.
Then
$$
\lim\limits_{t \to \infty} g\biggl[-\il{t_0}{t}\biggl(C(\tau) + A^*(\tau) B^{-1}(\tau) A(\tau)\biggr) d \tau - B^{-1}(t) A(t)\biggr] =\lim\limits_{t \to+\infty}g[-A_0] \ne \infty.
$$
Therefore, for this particular case Theorem 1.1 is not applicable to the system (1.1).  It is not difficult to verify that
$$
tr C(t) = - tr (A_0^* A_0) < 0, \ph t \ge t_0.
$$
Then since $A_0 + A^*_0 = 0$ using Theorem 3.2 to the system (1.1) we conclude that for this particular case the system (1.1) is oscillatory.}

Assume $B(t) \ge 0, \ph t \ge t_0$ and let $\sqrt{B_1(t)}, \ph t \ge t_0$ be absolutely continuous. Consider the linear matrix equation
$$
\sqrt{B_1(t)} X (A_1(t)\sqrt{B_1(t)} - \sqrt{B_1(t)}') = A_1(t)\sqrt{B_1(t)} - \sqrt{B_1(t)}', \ph t \ge t_0. \eqno (3.12)
$$
This equation has always a solution when $B_1(t) > 0, \ph t \ge t_0$. But it can have also a solution when $B(t)$ is not invertible for all (for some) $t \ge t_0$ (see [3,11]). In the general case Eq. (3.12) has a solution if and only if the equations
$$
\sqrt{B_1(t)}Y = A_1(t)\sqrt{B_1(t)} - \sqrt{B_1(t)}', \ph t \ge t_0, \phantom{aaaaaaaaaaaaaaaaaa}
$$
$$
 \phantom{aaaaaaaaaaaaaaaaaa} Z  A_1(t)\sqrt{B_1(t)} - \sqrt{B_1(t)}' =  A_1(t)\sqrt{B_1(t)} - \sqrt{B_1(t)}', \ph t\ge t_0
$$
have solutions (see [3], p. 23). Hence, Eq. (3.12) has a solution if and only if \linebreak $rank \sqrt{B_1(t)} = rank (\sqrt{B_1(t)} |  A_1(t)\sqrt{B_1(t)} - \sqrt{B_1(t)}'), \ph t \ge t_0.$

Let $F(t)$ be a solution of Eq. (3.12). We set:
$$
A_F(t) \equiv F(t)(A_1(t) \sqrt{B_1(t)} - \sqrt{B_1(t)}'), \phh J_2(t)\equiv -\frac{1}{2}tr (A_F(t) + A^*_F(t)) -
$$
$$
-\il{t_0}{t} tr \biggl[A_F(\tau) A^*_F(\tau)   + B(\tau) C(\tau)\biggr]d \tau +
\frac{1}{n}\il{t_0}{t}\Bigl[tr\Bigl(\frac{ A_F(\tau) - A^*_F(\tau)}{2i}\Bigr)\Bigr]^2 d\tau, \phh t \ge t_0.
$$

{\bf Theorem 3.3.} {\it Let $\sqrt{p(t)B(t)}$ be absolutely continuous on $[t_0,\infty)$ and  $F(t)$ be a solution of Eq. (3.12) such that the functions  $tr (A_F(t) A^*_F(t) + B(t) C(t)), \ph [tr( A_F(t) - A^*_F(t))]^2, \ph t \ge t_0$ are locally integrable on $[t_0,\infty)$. If
$$
\lim\limits_{t\to+\infty} J_2(t) = \infty
$$
then the system (1.1) is oscillatory.
}

Proof. Suppose the system (1.1) is not oscillatory. Then by (2.15) Eq. (2.14) has a solution $Y(t)$ on $[t_1,\infty)$ for some $t_1 \ge t_0$. Hence, by (2.16)
$$
Y'(t) + Y(t)B_1(t)Y(t) + A_1^*(t) Y(t) + Y(t) A_1(t) - C_1(t) = 0, \phh t \ge t_1.
$$
Multiply both sides of this equality at left and at right by $\sqrt{B_1(t)}$. Taking into account the equality
$$
(\sqrt{B_1(t)} Y(t) \sqrt{B_1(t)})' = \sqrt{B_1(t)}' Y(t) \sqrt{B_1(t)} + \sqrt{B_1(t)} Y'(t) \sqrt{B_1(t)} + \sqrt{B_1(t)} Y(t) \sqrt{B_1(t)}',
$$
$t \ge t_1$ we obtain
$$
(\sqrt{B_1(t)} Y(t) \sqrt{B_1(t)})' + (\sqrt{B_1(t)} Y(t) \sqrt{B(t)})^2 + (\sqrt{B_1(t)} A_1^*(t) - \sqrt{B_1(t)}') Y(t) \sqrt{B_1(t)} +\phantom{aaaaaaa}
$$
$$
\phantom{aaaaaaa}+ \sqrt{B_1(t)} Y(t) (A_1(t)\sqrt{B_1(t)} - \sqrt{B_1(t)}') - \sqrt{B_1(t)} C_1(t)\sqrt{B_1(t)}, \phh t \ge t_1.  \eqno (3.13)
$$
Since $F(t)$ is a solution of Eq. (3.12) we have $\sqrt{B_1(t)} A_1^*(t) - \sqrt{B_1(t)}' = (\sqrt{B_1(t)} A_1^*(t) - \sqrt{B_1(t)}') F^*(t) \sqrt{B_1(t)} = A^*_F(t)\sqrt{B_1(t)}= A^*_F(t) \sqrt{B_1(t)}, \ph A_1(t)\sqrt{B_1(t)} - \sqrt{B_1(t)}' = \linebreak A_F(t) \sqrt{B(t)},\ph t \ge t_1$. From here and  (3.13) it follows
$$
tr \{V(t)\}' + tr \{V^2(t) + A^*_F(t) V(t) + V(t) A_F(t) - \sqrt{B(t)} C(t) \sqrt{B(t)}\} =0, \ph t \ge t_1, \eqno (3.14)
$$
 where $V(t) \equiv \sqrt{B_1(t)} Y(t) \sqrt{B_1(t)}, \ph t \ge t_1.$
By Lemma 2.3
$$
tr[\sqrt{B(t)} C(t) \sqrt{B(t)}] = tr[B(t) C(t)], \phh t \ge t_0.
$$
Then if we substitute $V(t) \equiv Z(t) - \frac{A_F(t) + A^*_F(t)}{2}, \ph t \ge t_1$ in (3.14) (except in the expression $tr\{V(t)\}'$)  and take into account the condition of local integrability of $tr [A_F(t) A^*_F(t) + B(t) C(t)]$ and $[tr (A_F(t) - A^*_F(t))]^2$ we can, as in the proof of Theorem 3.1, to derive the inequality
$$
[tr Z(t)] + \il{t_1}{t}[tr Z(\tau)]^2 + J_2(t) + c_1 \le 0,
$$
where $c_1 \equiv - tr V(t_1) + \il{t_0}{t_1} tr [A_F(t) A^*_F(t) + B(t) C(t)] - \frac{1}{4 n} \il{t_0}{t_1} [tr (A_F(t) - A^*_F(t))]^2$  is a constant. Further as in the proof of Theorem 3.1. The theorem is proved.

{\bf Example 3.3.} {\it Assume $B(t) \equiv \frac{1}{p(t)}\begin{pmatrix} I_m & \theta_{12}\\ \theta_{21} & \theta_{22}\end{pmatrix}, \ph A(t)\equiv \begin{pmatrix} \theta_{11} & A_1(t)\\ \theta_{21} & A_2(t)\end{pmatrix} - \frac{1}{2}\bigl(\frac{p'(t)}{p(t)} + \linebreak +\mu(t)\bigr)I, \ph C(t)\equiv p(t) I, \ph rank A_2(t) \not\equiv 0, \ph t \ge t_0,$  where $I_m$ is an identity matrices of dimensions $m\times m$ ($m < n$), $\theta_{11}, \ph \theta_{12}, \ph \theta_{21}$ and $\theta_{22}$ are null matrices of dimensions $m\times m, \ph (n-m)\times m, \ph m\times (n-m)$ and $(m-m)\times (n-m)$ respectively.
Obviously for this case $F(t)\equiv 0, \ph t \ge t_0$ is a solution for Eq. (3.12). Then $A_F(t)\equiv 0, \phh t \ge t_0,, \ph J_2(t) =  (t - t_0) m \to +\infty$ for $t \to +\infty.$. By Theorem 3.3 it follows from here that in this particular case the system (1.1) is oscillatory.
}

\vskip 20 pt

\centerline{ \bf References}

\vskip 20pt

\noindent
1. K. I. Al - Dosary, H. Kh. Abdullah and D. Husein. Short note on oscillation of matrix \linebreak \phantom{a} Hamiltonian systems. Yokohama Math. J., vol. 50, 2003.

\noindent
2. Sh. Chen, Z. Zheng, Oscillation criteria of Yan type for linear Hamiltonian systems, \linebreak \phantom{a} Comput.  Math. with Appli., 46 (2003), 855--862.

\noindent
3.  F. R. Gantmacher, Theory of Matrix. Second Edition (in Russian). Moskow,,\linebreak \phantom{a} ''Nauka'', 1966.

\noindent
4. G. A. Grigorian, Oscillation criteria for linear matrix Hamiltonian systems. \linebreak \phantom{aa} Proc. Amer. Math. Sci, Vol. 148, Num. 8 ,2020, pp. 3407 - 3415.

\noindent
5.  G. A. Grigorian.   Interval oscillation criteria for linear matrix Hamiltonian systems,\linebreak \phantom{a}  Rocky Mount. J. Math.,  vol. 50 (2020), No. 6, 2047–2057

\noindent
6.   G. A. Grigorian. Oscillatory criteria for the systems of two first - order Linear \linebreak \phantom{a} ordinary differential equations. Rocky Mount. J. Math., vol. 47, Num. 5,
 2017, \linebreak \phantom{a}  pp. 1497 - 1524

\noindent
7. G. A. Grigorian,  On two comparison tests for second-order linear  ordinary\linebreak \phantom{aa} differential equations (Russian) Differ. Uravn. 47 (2011), no. 9, 1225 - 1240; trans-\linebreak \phantom{aa} lation in Differ. Equ. 47 (2011), no. 9 1237 - 1252, 34C10.

\noindent
8.  G. A. Grigorian, Two comparison criteria for scalar Riccati equations and their \linebreak \phantom{aa} applications.
Izv. Vyssh. Uchebn. Zaved. Mat., 2012, Number 11, 20–35.

\noindent
9.  G. A. Grigorian, Oscillatory and Non Oscillatory criteria for the systems of two \linebreak \phantom{aa}   linear first order two by two dimensional matrix ordinary differential equations. \linebreak \phantom{aa}   Arch.  Math., Tomus 54 (2018), PP. 189 - 203.

\noindent
10.  I. S. Kumary and S. Umamaheswaram, Oscillation criteria for linear matrix \linebreak \phantom{a} Hamiltonian systems, J. Differential Equ., 165, 174--198 (2000).

\noindent
11.  A. Kurosh, Higher Algebra, Moskow Mir Publisher (English translation) 1980, \linebreak \phantom{a}428 pages.

\noindent
12.  L. Li, F. Meng and Z. Zheng, Oscillation results related to integral averaging technique\linebreak \phantom{a} for linear Hamiltonian systems, Dynamic Systems  Appli. 18 (2009), \ph \linebreak \phantom{a} pp. 725--736.

\noindent
13.  Y. G. Sun, New oscillation criteria for linear matrix Hamiltonian systems. J. Math. \linebreak \phantom{a} Anal. Appl., 279 (2003) 651--658.

\noindent
14. C. A. Swanson. Comparison and oscillation theory of linear differential equations.   \linebreak \phantom{a} Academic press. New York and London, 1968.

\noindent
15. Q. Yang, R. Mathsen and S. Zhu, Oscillation theorems for self-adjoint matrix \linebreak \phantom{a}   Hamiltonian
 systems. J. Diff. Equ., 19 (2003), pp. 306--329.

\noindent
16.  Z. Zheng, Linear transformation and oscillation criteria for Hamiltonian systems. \linebreak \phantom{a} J. Math. Anal. Appl., 332 (2007) 236--245.

\noindent
17. Z. Zheng and S. Zhu, Hartman type oscillatory criteria for linear matrix Hamiltonian  \linebreak \phantom{a} systems. Dynamic  Systems  Appli., 17 (2008), pp. 85--96.

\end{document}